\documentclass[12pt]{article}

\usepackage[cp1251]{inputenc}
\usepackage[english]{babel}
\usepackage{amsmath,amsfonts,amssymb}
\usepackage{graphics}
\author{Kochkarev B.S.}
\title{Absolutely symmetric trees and complexity of natural number }
\begin{document}
\maketitle
Abstract. We consider the rooted trees which not have isomorphic
representation and  introduce a conception of complexity a natural number also.
The connection between quantity such trees with $n$ edges and a complexity of natural number $n$ is established. The recurrent ratio for complexity of a natural number is founded. An expression for calculation of difference complexities of two adjacent natural numbers is constructed. It is proved that this difference equal 1 if and only if a natural number is simple. From proved theorems it follows corollaries.

Almost every book on graph theory contains some parts devoted to trees (see, e.g., \cite{1,2,3,4}).
 The concept of a tree has been entered for the first time by Kirchhoff \cite{5} in connection with research of electric chains. Later this concept also was independent is entered by Cayley \cite{6} and it had been received the first basic results in the theory of trees. Trees have wide appendices in various area of a science. The analysis of publications [1,2,3,4,5,6,7,8,9,10,11,12,13,14,15,16] allows to assume that a class of rooted trees which not have isomorphic trees distinct from them in the scientific literature is disregarded. The connection between quantity such trees with $n$ edges and a complexity of natural number $n$ is established in this article.

Definition 1. A flat geometrical realization of a rooted tree is called absolutely symmetric rooted tree (a.s.r.t.) if it not have isomorphic representative distinct from them.

Evidently at $n=1,2$ all flat geometrical realization of rooted trees are a.s.r.t.. Since $n=3$ among realization of trees with $n$ edges there are isomorphic trees.

It is possible to make other inductive definition a.s.r.t. equivalent given above.

Induction basis. The edge with the allocated end $a$ (fig.1) is a.s.r.t. with a root $a$.

Induction step. Let $A$ - be a a.s.r.t. with a root $a$ (fig.2). Then the tree $C$ (fig.3a) obtained from $A$ by "connecting" a new edge to the root $a$, become a a.s.r.t. with the root $c$, where $c$ is the free end of the connected edge. Further, the tree $D$ (fig.3b), obtained from $k, k>1, A$ by joining their roots, become a a.s.r.t. with the root $a$.

Let $A=(a_{1},a_{2},\ldots,a_{k})$ a vector with $k$ natural components $a_{i}\geq 1$. For $A$ we will define functional $f(A)$ recurrently: $f(a_{1})=a_{1},f(a_{1},a_{2},\ldots,a_{i+1})=(f(a_{1},a_{2},\ldots,a_{i})+1)a_{i+1}$.

Definition 2. The quantity $T(n)$ of various vectors $A$ with $f(A)=n$ is called a complexity of number $n$.

Theorem 1. The quantity of a.s.r.t. with $n$ edges equal $T(n)$.

 Proof. We show at first to any vector $A=(a_{1},a_{2},\ldots,a_{k}), a_{i}\geq1$, it is
 possible to put in conformity a.s.r.t. with $f(A)$ edges. We prove it an induction on $k$.
 For $k=1$  $A=(a_{1}), f(A)= a_{1}$. We put to vector $A$ in conformity in this
 case a rooted tree $A^{'}_{1}$ (fig.4) which is a.s.r.t. with $a_{1}$ edges and with a root
 $a$. We assume now, that for any $l\leq i$ to vector $A=(a_{1},a_{2},\ldots,a_{l})$
 corresponds a.s.r.t. $A^{'}_{l}$ with number of edges equal $f(A)$  and with a root $a$.
 We will show how, that to vector $A=(a_{1},a_{2},\ldots,a_{i+1})$ can put in conformity a.s.r.t.
 $A^{'}_{i+1}$ with number of edges $f(A)$.Really, under the assumption, to vector
 $A=(a_{1},a_{2},\ldots,a_{i})$ corresponds a.s.r.t. $A^{'}_{i}$ (fig.5a) with number of
 edges equal $f(A)$ and with a root $a$. Then to vector $A=(a_{1},a_{2},\ldots,a_{i},a_{i+1})$
 we will put in conformity rooted tree $A^{'}_{i+1}$ with a root $c$ (fig.5b) obtained from
 $a_{i+1}$ trees $A^{''}_{i}$ (fig.5c) by joining their roots $c$. Evidently, $A^{'}_{i+1}$
 is a.s.r.t. with number of edges equal $f(A)=(f(a_{1},a_{2},\ldots,a_{i})+1)a_{i+1}$, as was to be shown.
 Thus to any vector $A=(a_{1},a_{2},\ldots,a_{k}), a_{i}\geq 1$ it is possible to put in
 conformity a.s.r.t. with $f(A)$ edges.

We show now, that to any a.s.r.t. with number of edges $n$ there correspond some vector $A=(a_{1},a_{2},\ldots,a_{k})$ wits $f(A)=n$. We prove it an induction on $n$. If $n=1$, a.s.r.t. with number of edges 1 (fig.1) correspond vector $A=(1)$ with $f(A)=1$ We assume now, that for any $l\leq k$ a.s.r.t. with number of edges $l$ there correspond some vector $A=(a_{1},a_{2},\ldots,a_{m})$ with $f(A)=l$. Let now it is given a.s.r.t. with number of edges $k+1$. To inductive definition a.s.r.t. two cases are possible: the tree looks like, presented on fig.6a or 6b, where $A$ and $B$ a.s.r.t. such, that vectors corresponding to them $A^{'}=(a_{1},a_{2},\ldots,a_{s}), B^{'}=(a^{'}_{1},a^{'}_{2},\ldots,a^{'}_{r})$ have $f(A^{'})=k,f(B^{'})=m$. Then a.s.r.t. $C$(fig.6a) and $D$(fig.6b) there will correspond vectors $A_{C}=(a_{1},a_{2},\ldots,a_{s},1),A_{D}=(a^{'}_{1},a^{'}_{2},\ldots,a^{'}_{r},t)$.
Evidently, $f(A_{C})=k+1,f(A_{D})=(f(B^{'})+1)t=k+1$. Thus, to any a.s.r.t. with number of edges $n$ there correspond some vector $A=(a_{1},a_{2},\ldots,a_{k})$ with $f(A)=n$.

Corollary 1. If $n$ is simple then all vectors $A=(a_{1},a_{2},\ldots,a_{k}),k>1$, for which $f(A)=n$ have $a_{k}=1$.

Corollary 2. $T(n)$ is the quantity of presentations of $n$ in look $a_{1}a_{2}\cdots a_{k}+a_{2}a_{3}\cdots a_{k}+ \cdots+a_{k}$.

Theorem 2. $T(n)=\sum_{r} T(r-1)$, where $r\geq 1$ divide $n,T(0)=1$.

Proof. We prove it an induction on $n$. $T(1)=1$ (fig.1). Really, $T(1)=\sum_{r}T(r-1)=T(0)=1$.We assume now, that the statement is true for all $l\leq k$, i.e. $T(l)=\sum_{r}T(r-1)$. Let now $l=k+1$. A.s.r.t. with $k+1$ edges can be one of two kinds (fig.6a) or(fig.6b). Obviously, the number of a.s.r.t a kind (fig.6a) will be $T(k)$, and the number of a.s.r.t. a kind (fig.6b) will be $\sum_{r'}T(r'-1)$, where $r'$ divide $k+1$ and $r'<k+1$. Thus, the general number of a.s.r.t. with $k+1$ edges will be $T(k+1)=T(k)+\sum_{r'} T(r'-1)=\sum_{r}T(r-1)$. From here for any $n$ it is had

$T(n)=\sum_{r}T(r-1)$ (1)

For example, if $n=12$ then $T(12)= T(0)+T(1)+T(2)+T(3)+T(5)+T(11)=40$.

Corollary 1. Number $n$ is simple if and only if $T(n)= 1+T(n-1)$

Corollary 2. If $n=p^{r}$, where $p$ is simple then $T(n)=\sum_{k=0}^{r}T(p^{k}-1)$.

Corollary 3. If $n$ is compound then $T(n)<1+T(n-1)+(2\lfloor\sqrt{n}\rfloor-2)
T\big(\frac{n}{\min\limits_{i}p_{i}}-1\big)$, where $p_{i}$ is simple divider of $n$.

Corollary 4. $T(n)\leq{n \choose \lfloor\frac{n}{2}\rfloor}$.

Proof. It statement is proved on induction. For $n=1$ it is true. We assume now that it is true for all
 $l\leq k$, i.e.
 $T(l)\leq{l\choose \lfloor\frac{l}{2}\rfloor}$.
 Then from corollary 2 we receive
$T(k+1)< 1+T(k)+(2\lfloor\sqrt{k+1}\rfloor -2)
T\big(\frac{k+1}{\min\limits_{i}p_{i}}-1\big)
\leq{k+1\choose\lfloor\frac{k+1}{2}\rfloor}$.
Thus, for any $n$ we have $T(n)\leq {n \choose \lfloor \frac{n}{2}\rfloor}$.

Let $T^\ast(n)=T(n)-T(n-1),n\geq1$. Evidently, for $n\geq2$ $T^\ast(n)$ represent a difference on complexity of two adjacent natural numbers, i.e. on how many complexity of natural number $n$ is more than complexity of $n-1$.

Theorem 3. $T^{\ast}(n)=\sum_{p_{i}}T(\frac{n}{p_{i}})-\sum_{p_{i},p_{j}}T(\frac{n}{p_{i}p_{j}})+\cdots+(-1)^{m-1}T(\frac{n}{p_{1}p_{2}\cdots p_{m}})$,
where $p_{1},p_{2},\ldots,p_{m}$ runs all in pairs various simple dividers of number $n\geq2$

Proof. According to (1) $T(n)=\sum_{r}T(r-1)=\sum_{r'}T(r'-1)+T(n-1)$ From here $T^\ast(n)=\sum_{r'}T(r'-1)$ It is obvious, if $n=p^{r_{1}}_{1}p^{r_{2}}_{2}\cdots p^{r_{m}}_{m}$, where $p_{1},p_{2},\ldots,p_{m}, m\geq2$, runs all in pairs various simple dividers of number $n$ then $\sum_{i=1}^{m}T(\frac{n}{p_{i}})\geq\sum_{r'}T(r'-1)$ as any divider $r'$ of number $n$ distinct from $n$ is a divider at least one of numbers $\frac{n}{p_{i}}, i=\overline{1,m}$. At numbers $\frac{n}{p_{i}}, i=\overline{1,m}$, can be the general dividers. There fore, using a method of inclusions and exeption [16] we receive $T^{\ast}(n)=\sum_{p_{i}}T(\frac{n}{p_{i}})-\sum_{p_{i},{p_[j}}T(\frac{n}{p_{i}p_{j}})+ \cdots +(-1^{m-1}T(\frac{n}{p_{1}p_{2}\cdots p_{m}})$.

Corollary. $T^{\ast}(n)=1$ if and only if $n$ is  simple.

Proof. Let $n$ is simple. Then $T^{\ast}=T(n)-T(n-1)=T(0)+T(n-1)-T(n-1)=T(0)=1$. Let now $T^{\ast}=1$. Then as $T^{\ast}=T(n)-T(n-1)=1$, $T(n)=1+T(n-1)$. But $T(n)=\sum_{r}T(r-1)=T(0)+T(n-1)+\sum_{r''}T(r''-1)$, where $r''>1$ divide $n$ and $r''\neq n$. From here $\sum_{r''}T(r''-1)=0$ Hence $n$ is simple.

From definition $T^{\ast}(n),n\geq1$, we have $T(n)=1+\sum_{k=1}^{n}T^\ast(k)$ The first 40 values of numbers $T^{\ast}(n), n=1,2,\ldots,40$, are that 0,1,1,2,1,4,1,5,3,7,1,13,1,12,
8,16,1,26,1,29,13,28,1,51,6,42,19,56,1,87,1,77,29,79,16,124,1,106,43,145.

References:

\newpage

\begin{figure}[h]
 \begin{center}
 \includegraphics {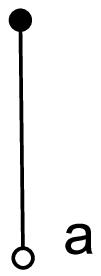}
 \caption{}
 \label{Ris1}
  \end{center}
\end{figure}

\begin{figure}[h]
 \begin{center}
 \includegraphics {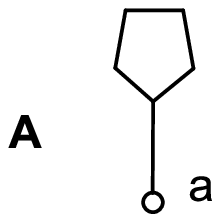}
 \caption{}
 \label{Ris2}
  \end{center}
\end{figure}

\begin{figure}[htb]
 \begin{center}
 \vbox{ \scalebox{0.22}{ \includegraphics{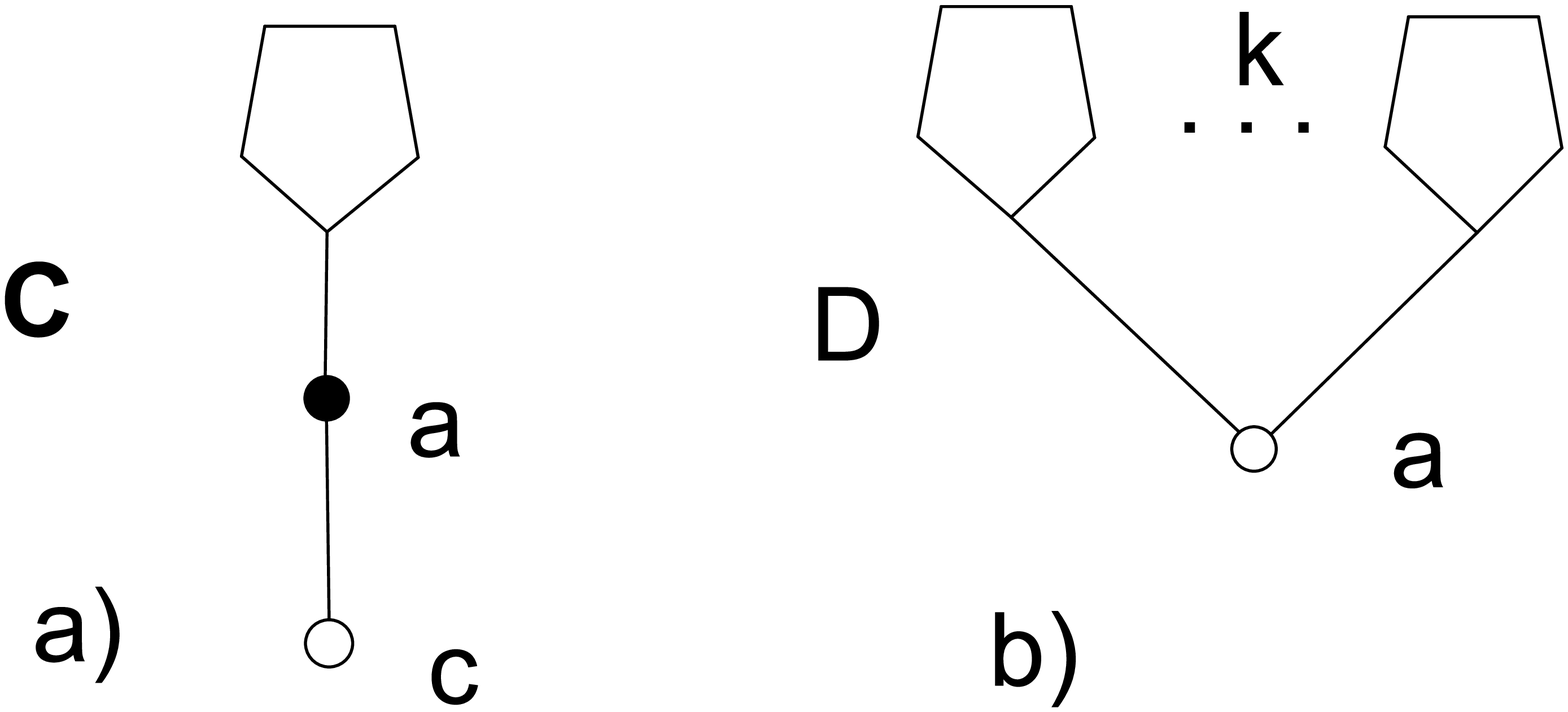}} }
 \caption{}
 \label{Ris3}
  \end{center}
 \end{figure}

\begin{figure}[htb]
 \begin{center}
 \vbox{ \scalebox{0.22}{ \includegraphics{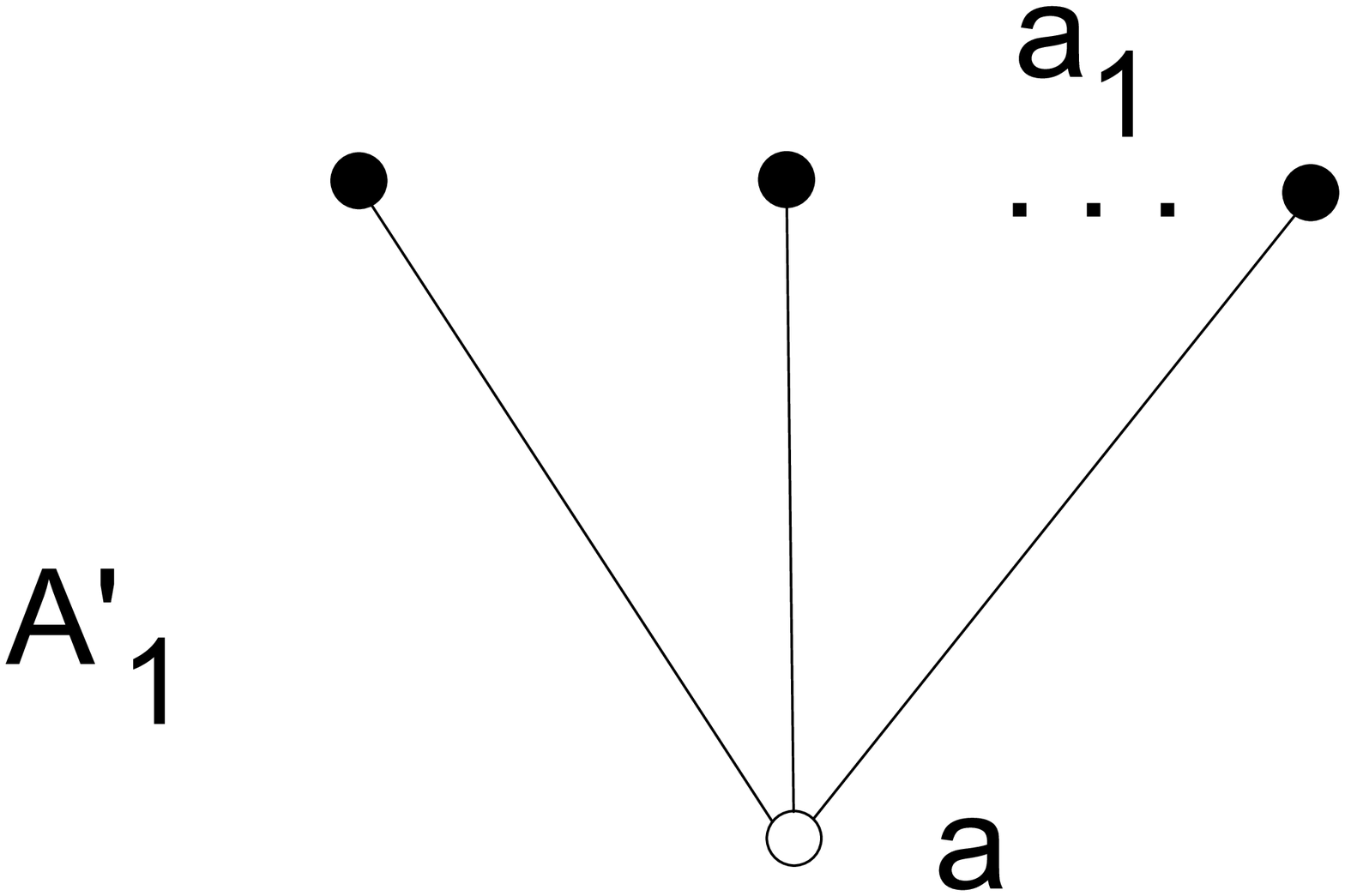}} }
 \caption{}
 \label{Ris4}
  \end{center}
 \end{figure}

\begin{figure}[htb]
 \begin{center}
 \vbox{ \scalebox{0.22}{ \includegraphics{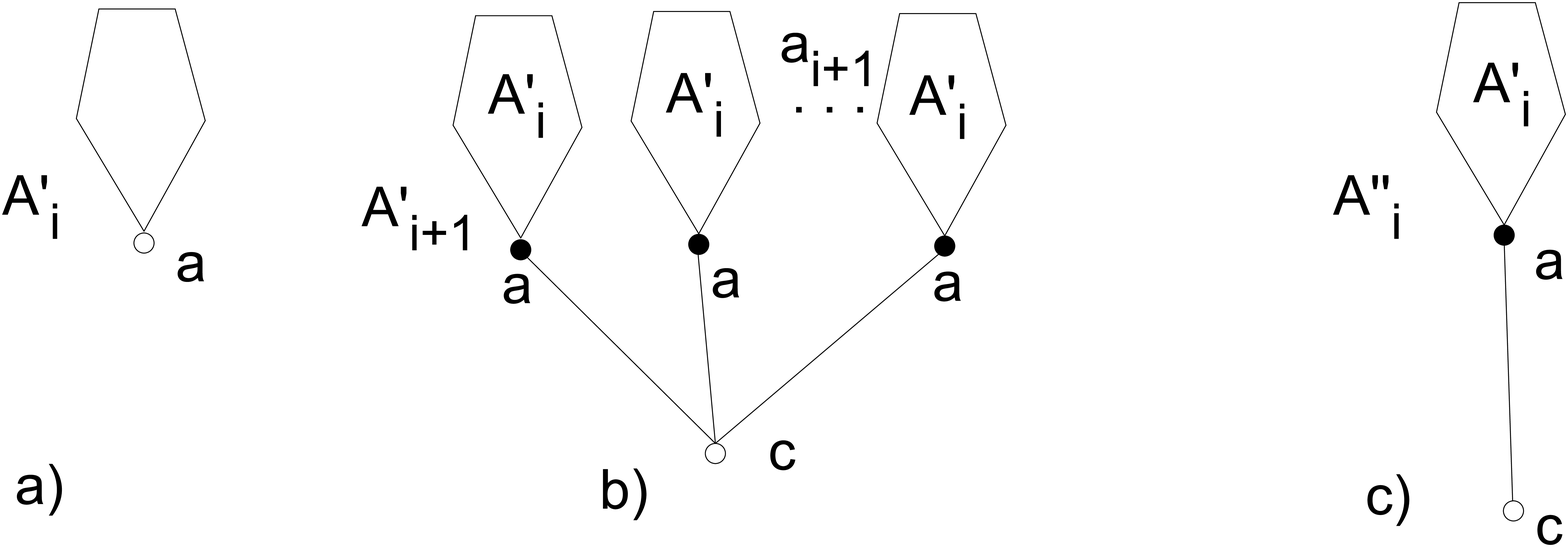}} }
 \caption{}
 \label{Ris5}
  \end{center}
 \end{figure}

\begin{figure}[htb]
 \begin{center}
 \vbox{ \scalebox{0.22}{ \includegraphics{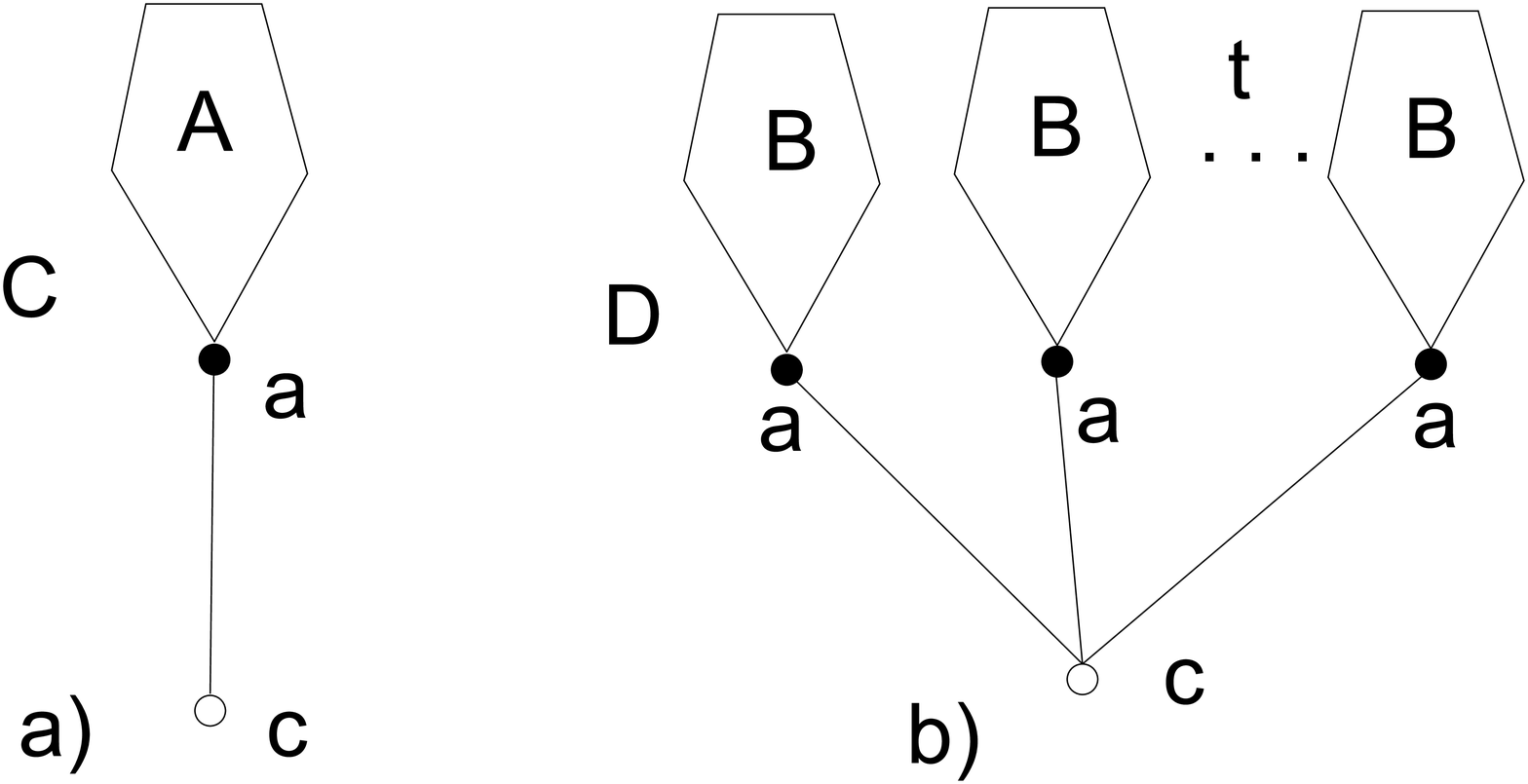}} }
 \caption{}
 \label{Ris6}
  \end{center}
 \end{figure}

\end{document}